\documentclass[12 pt]{article}
\usepackage{graphicx} 
\usepackage{graphicx}
\usepackage{amsmath,amsthm,amssymb,enumerate}
\usepackage{euscript,mathrsfs}
\usepackage{color}
\usepackage{dsfont}
\usepackage[left=2cm,right=2cm,top=3.5cm,bottom=3.5cm]{geometry}
\usepackage{color}
\usepackage[framemethod=tikz]{mdframed}
\allowdisplaybreaks

\usepackage{soul}
\usepackage{esint}
\catcode`\@=11 \@addtoreset{equation}{section}

\catcode`\@=12

\allowdisplaybreaks

\newtheorem{Theorem}{Theorem}[section]
\newtheorem{Proposition}[Theorem]{Proposition}
\newtheorem{Lemma}[Theorem]{Lemma}
\newtheorem{Corollary}[Theorem]{Corollary}

\theoremstyle{definition}
\newtheorem{Definition}[Theorem]{Definition}

\newtheorem{Remark}[Theorem]{Remark}

\newcommand{\bTheorem}[1]{
\begin{Theorem} \label{T#1} }
\newcommand{\eT}{\end{Theorem}}

\newcommand{\bProposition}[1]{
\begin{Proposition} \label{P#1}}
\newcommand{\eP}{\end{Proposition}}

\newcommand{\bLemma}[1]{
\begin{Lemma} \label{L#1} }
\newcommand{\eL}{\end{Lemma}}

\newcommand{\bCorollary}[1]{
\begin{Corollary} \label{C#1} }
\newcommand{\eC}{\end{Corollary}}

\newcommand{\bRemark}[1]{
\begin{Remark} \label{R#1} }
\newcommand{\eR}{\end{Remark}}

\newcommand{\bDefinition}[1]{
\begin{Definition} \label{D#1} }
\newcommand{\eD}{\end{Definition}}

\newcommand{\bfphi}{\boldsymbol{\varphi}}

\newcommand{\bFormula}[1]{
\begin{equation} \label{#1}}
\newcommand{\eF}{\end{equation}}

\newcommand{\Ov}[1]{\overline{#1}}

\newcommand{\vr}{\varrho}

\newcommand{\tvr}{\tilde \vr}

\newcommand{\vt}{\vartheta}

\newcommand{\vm}{\vc{m}}

\newcommand{\vc}[1]{{\bf #1}}

\newcommand{\Div}{{\rm div}_x}
\newcommand{\Grad}{\nabla_x}

\newcommand{\dx}{\,{\rm d} {x}}

\newcommand{\dt}{\,{\rm d} t }

\newcommand{\vU}{\vc{U}}

\newcommand{\intO}[1]{\int_{\Omega} #1 \ \dx}

\newcommand{\D}{{\rm d}}

\newcommand{\ep}{\varepsilon}

\newcommand{\br}{ \nonumber \\ }

\def\softd{{\leavevmode\setbox1=\hbox{d}%
          \hbox to 1.05\wd1{d\kern-0.4ex{\char039}\hss}}}
\definecolor{Cgrey}{rgb}{0.85,0.85,0.85}
\definecolor{Cblue}{rgb}{0.50,0.85,0.85}
\definecolor{Cred}{rgb}{1,0,0}
\definecolor{fancy}{rgb}{0.10,0.85,0.10}

\newcommand\Cbox[2]{%
    \newbox\contentbox%
    \newbox\bkgdbox%
    \setbox\contentbox\hbox to \hsize{%
        \vtop{
            \kern\columnsep
            \hbox to \hsize{%
                \kern\columnsep%
                \advance\hsize by -2\columnsep%
                \setlength{\textwidth}{\hsize}%
                \vbox{
                    \parskip=\baselineskip
                    \parindent=0bp
                    #2
                }%
                \kern\columnsep%
            }%
            \kern\columnsep%
        }%
    }%
    \setbox\bkgdbox\vbox{
        \color{#1}
        \hrule width  \wd\contentbox %
               height \ht\contentbox %
               depth  \dp\contentbox
        \color{black}
    }%
    \wd\bkgdbox=0bp%
    \vbox{\hbox to \hsize{\box\bkgdbox\box\contentbox}}%
    \vskip\baselineskip%
}

\mdfdefinestyle{MyFrame}{%
	linecolor=black,
	outerlinewidth=1pt,
	roundcorner=5pt,
	innertopmargin=\baselineskip,
	innerbottommargin=\baselineskip,
	innerrightmargin=10pt,
	innerleftmargin=10pt,
	backgroundcolor=white!20!white}

\date{}

\begin{document}


\title{\bf On ``discrete'' solutions of the Euler system of gas dynamics}

\author{Anna Abbatiello\footnote{\texttt{anna.abbatiello@unicampania.it}} 
\and Eduard Feireisl\footnote{\texttt{feireisl@math.cas.cz}}}

\date{}

\maketitle

\medskip

\centerline{University of Campania ``L.~Vanvitelli", Department of Mathematics and Physics}

\centerline{Viale A.~Lincoln 5, 81100 Caserta, Italy}

\medskip

\centerline{Institute of Mathematics of the Academy of Sciences of the Czech Republic,}

\centerline{\v Zitn\' a 25, CZ-115 67 Praha 1, Czech Republic}

\maketitle

\begin{abstract}
	
The method of Convex Integration has revealed a number of 
rather disturbing facts concerning well-posedness of the Euler system of gas dynamics. In particular, there is a dense set of 
``wild'' initial data, for which the problem admits 
 infinitely many physically admissible (entropy) weak solutions.
 We identify the class of initial data enjoying the following properties:
 \begin{itemize}
 \item they give rise to a family of weak solutions with increasing 
 entropy profiles;
 \item the solutions are ``discrete'', meaning they attain only a finite number 
 of constant states;
 \item the solutions reach a prescribed terminal entropy profile when time goes to infinity.
 \end{itemize}
 
\end{abstract}

\bigskip


\bigskip

\section{Introduction}

The Euler system of gas dynamics is 
one of the simplest models in continuum mechanics governing
the time evolution of a compressible inviscid liquid like a gas. 
It can be formally seen as the zero dissipation limit of problems describing 
real fluids featuring viscosity and heat diffusion. Despite its apparent simplicity, the structure of solutions of the Euler system is still quite poorly understood, in particular in the multidimensional setting. Smooth solutions develop singularities in a finite time, while the problem is ill--posed in the larger class of weak solutions even if some 
physically relevant admissibility criteria are imposed. 

\subsection{Euler system}

The time evolution of the phase variables, the mass density $\vr= \vr(t,x)$, the momentum $\vm = \vm(t,x)$, and the energy $E= E(t,x)$ 
is described by means of a system of field equations:

\medskip

\noindent {\bf Mass conservation:} 
\begin{equation} \label{i1}
	\partial_t \vr + \Div \vm = 0; 
\end{equation}	

\noindent {\bf Momentum balance:}
\begin{equation} \label{i2}
\partial_t \vm + \Div \left( \frac{\vm \otimes \vm}{\vr} \right) + \Grad p = 0;
\end{equation}	

\noindent {\bf Energy balance:}
\begin{equation} \label{i3}
\partial_t E + \Div \left[ \left(E + p \right) \frac{\vm}{\vr} \right] = 0,
\end{equation}
where $p$ is the pressure. In addition, we introduce the absolute temperature 
$\vt$, the internal energy $e$, and the entropy $s$ through Gibbs' relation
\begin{equation} \label{i4}
\vt Ds = De + p D \left( \frac{1}{\vr} \right).
\end{equation}
It is convenient to pass to new state variables $(\vr, s, \vm)$, and write 
\begin{align} 
p = p(\vr, s), \
E 	= E(\vr, s, \vm) = \frac{1}{2} \frac{|\vm|^2}{\vr} + \vr e (\vr, s).
\label{i5}
\end{align}	
Finally, the Second Law of Thermodynamics is enforced through {\bf entropy inequality:} 
\begin{equation} \label{i6}
\partial_t (\vr s) + \Div (s \vm) \geq 0.
\end{equation}

We suppose the fluid is confined to a bounded domain $\Omega \subset R^d$, $d=2,3$, and we consider the time evolution of the system for $t \in [0,T)$, $T \leq \infty$. Furthermore, we assume the boundary $\partial \Omega$ is impermeable, meaning 
\begin{equation} \label{i7}
	\vm \cdot \vc{n}|_{\partial \Omega} = 0, 
\end{equation}
where $\vc{n}$ denotes the outer normal vector. The original state of the system is determined by the initial conditions
\begin{equation} \label{i8}
\vr(0, \cdot) = \vr_0,\ \vm(0, \cdot) = \vm_0, \ 
s(0, \cdot) = s_0 \ \Rightarrow \ 
E(0, \cdot) = \frac{1}{2} \frac{|\vm_0|^2}{\vr_0} + \vr_0 e(\vr_0, s_0).
\end{equation}
Here and hereafter, we tacitly assume $\vr_0 > 0$ to avoid problems with vacuum zones. 

The present state-of-the-art of the mathematical theory of the Euler system includes the following facts: 
\begin{itemize}
	\item The Euler system is well--posed, locally in time, in the class of sufficiently regular initial data, see e.g.~Benzoni-Gavage and Serre \cite[Part IV] {BenSer} or Schochet \cite[Theorem 1]{SCHO1}.
	
	\item Regular solutions may develop singularities in the form of shock waves appearing in a finite time, see e.g.~the monograph of Smoller 
	\cite{SMO}.
	
	\item Regular solutions may develop implosions (blow-up type singularities) 
	in a finite time, see Chen et al. \cite{ChCiShVi26}, \cite{ChShVi}.
	
	\item  The theory of convex integration
	developed in the context of fluid mechanics in the seminal work of De Lellis and Sz{\'e}kelyhidi \cite{DelSze13} and \cite{DelSze3},
	revealed several rather surprising facts concerning well--posedness of the Euler system even in the class of \emph{entropy admissible} 
	weak solutions. The ill--posedness of the Riemann problem in the 
	physically relevant 2 - 3--d geometry was shown by Klingenberg et al.~\cite[Theorem 1.1]{KlKrMaMa}, see also Al Baba et al.~\cite{ABKlKrMaMa}. 
	
	\item There is a large class of so-called ``wild'' initial data that give rise to infinitely many entropy admissible weak solutions, see Feireisl et al.   \cite{FeKlKrMa}.
	
	\item The wild data are dense in the $L^p-$topology of the phase space, see Chiodaroli, Feireisl \cite{ChFe2023}.
\end{itemize}			

\subsection{Discrete solutions of the Euler system}

Luo, Xie, and Xin \cite{LuXiXi} showed that the simplified isentropic (barotropic) 
Euler system admits weak solutions that may attain only a finite number of constant states. We call these solutions \emph{discrete}. Our aim is to extend this result to the full Euler system of gas dynamics \eqref{i1}--\eqref{i3}, where, in addition, the method will generate a family of solutions with increasing entropy production rate. More specifically, our goal is to show that 
there is a set $I_{\vr,s}$ of initial densities and entropies $(\vr_0, s_0)$ dense in the 
$L^p-$topology enjoying the following properties:
\begin{itemize}
\item For each $(\vr_0, s_0) \in I_{\vr,s}$, there exists $\vm_0 \in L^\infty(\Omega; R^d)$ such that the Euler system \eqref{i1}--\eqref{i8} 
admits an infinite family of weak solutions $(\vr^\lambda, \vm^\lambda, s^\lambda)_{\lambda \geq 0}$ in $(0,T) \times \Omega$ emanating from the initial data $(\vr_0, \vm_0, s_0)$.
\item For each $\lambda \geq 0$, there exists a finite integer $K$ such that 
\begin{equation} \label{i9}
(\vr^\lambda, \vm^\lambda, s^\lambda) (t,x) \in \left\{ (\vr_k, \vm_k, s_k) \Big| 
k = 1, \dots, K \right\} \ \mbox{for a.a.}\ (t,x) \in (0,T) \times \Omega.
\end{equation}
for some constant vectors $(\vr_k, \vm_k, s_k) \in (0, \infty) \times R^d \times R$.

\item 
\begin{align} 
\intO{ \vr^{\lambda_1} s^{\lambda_1} (t, \cdot) } &\leq 
\intO{ \vr^{\lambda_2} s^{\lambda_2} (t, \cdot) } \ \mbox{for a.a.}\ t \in (0,T),\br 
\intO{ \vr^{\lambda_1} s^{\lambda_1}  } &\ne \intO{ \vr^{\lambda_2} s_{\lambda^2}  } \ \mbox{in}\ (0,T) 
\label{i10}
\end{align}
whenever $\lambda_1 < \lambda_2$.

\end{itemize}	

Property \eqref{i9} means that the solutions are discrete - they attain only a finite number of states. Property \eqref{i10} says that 
these solutions fail to be maximal neither in DiPerna's nor in Dafermos' sense 
as specified in Feireisl, Luk\' a\v cov\' a, and Yu \cite{FeiLukYu}.

Alternatively, instead of prescribing the initial distribution of the density 
and the entropy, we may fix the initial density profile $\vr_0$, together with the initial kinetic energy 
\[
E_{0, {\rm kin}} = \frac{1}{2} \frac{|\vm_0|^2}{\vr_0}.
\]

Our final objective is to show that the same set of initial data give rise 
to a family of solutions attaining a given terminal entropy profile $\widetilde{s}$, meaning
\[
\vr s (t, \cdot) \nearrow \widetilde{s} \ \mbox{as}\ t \to T.
\]

The paper is organized as follows. In Section \ref{m}, we introduce the concept of admissible weak solution to the Euler system and state our main results. 
Next, we reformulate the problem and reduce the task to the investigation 
of the incompressible Euler system, see Section \ref{r}. In Section \ref{d}, 
we adapt the result of Luo, Xie, Xin \cite{LuXiXi} to the transformed problem. 
We complete the proof of the main results by constructing solutions with the increasing entropy production rate, and solutions with a given terminal entropy in Section \ref{e}. The paper is concluded by Section \ref{C}, where
the impact of the main results on the long time behaviour of the system is discussed, together with their admissibility in a larger class of 
dissipative (measure--valued) solutions, see e.g. Feireisl et al. \cite{FeLMMiSh}.

\section{Weak formulation and the main results}
\label{m}

Throughout the whole paper, we consider solutions with positive density $\vr$. In particular, all terms including $\vr^{-1}$ are well defined. For the sake of simplicity, we suppose the gas is \emph{polytropic}, meaning 
\begin{equation} \label{m1}
p = (\gamma - 1) \vr e,\ \mbox{with the adiabatic constant}\ 1 < \gamma \leq \frac{5}{3}.
\end{equation}	

\begin{Definition}[\bf Admissible weak solution] \label{mD1}
Let the initial data 
\begin{equation} \label{m2}
\vr_0\ \mbox{measurable},\ 0 < \underline{\vr} \leq \vr_0(x) \leq 
\Ov{\vr}\ \mbox{for a.a.}\ x \in \Omega,\ \vm_0 \in L^\infty(\Omega; R^d),\ s_0 \in L^\infty(\Omega; R^d)
\end{equation}	
be given. We say that a trio $(\vr, \vm, s) \in L^\infty((0,T) \times \Omega; 
R^{d+2})$ is \emph{admissible weak solution} to the Euler system \eqref{i1}--\eqref{i8} if the following holds:
\begin{itemize}
\item {\bf Absence of vacuum zones:}
\begin{equation} \label{m3}
	{\rm ess} \inf_{(0,T) \times \Omega} \vr > 0.
\end{equation}
\item {\bf Mass conservation:}
\begin{equation} \label{m4}
\int_0^T \intO{ \Big[ \vr \partial_t \varphi + \vm \cdot \Grad \varphi \Big] } \dt = 
- \intO{ \vr_0 \varphi(0,\cdot)} 
\end{equation}	
for any $\varphi \in C^1_c([0,T) \times \Ov{\Omega})$.
\item {\bf Momentum balance:}
\begin{equation} \label{m5}
\int_0^T \intO{ \left[ \vm \cdot \partial_t \bfphi + 
	\left(	\frac{\vm \otimes \vm }{\vr} \right) : \Grad \bfphi + 
	p(\vr,s) \Div \bfphi \right] } \dt = - \intO{ \vm_0 \cdot \bfphi(0,\cdot) } 
\end{equation}
for any $\bfphi \in C^1_c([0,T) \times \Ov{\Omega}; R^d)$, $\bfphi \cdot \vc{n}|_{\partial \Omega} = 0$.
\item {\bf Energy balance:}
\begin{align} 
\int_0^T &\intO{ \left[ E(\vr, s, \vm) \partial_t \varphi + 
\Big( E(\vr, s, \vm) + p(\vr,s) \Big) \frac{\vm}{\vr} \cdot \Grad \varphi \right]	} \dt \br & = - \intO{ 	E(\vr_0, s_0, \vm_0) \varphi (0, \cdot) } 
\label{m6}
\end{align}				
for any $\varphi \in C^1_c([0,T) \times \Ov{\Omega})$.
\item {\bf Entropy inequality:}
\begin{equation} \label{m7}
\int_0^T \intO{ 
\Big[ \vr s \partial_t \varphi + s \vm \cdot \Grad \varphi \Big] } \dt 
\leq - \intO{ \vr_0 s_0 \varphi (0, \cdot) } 
\end{equation} 
for any $\varphi \in C^1_c([0,T) \times \Ov{\Omega})$, $\varphi \geq 0$.		
	
\end{itemize}

\end{Definition}

\begin{Remark} \label{mR2}
	
Since no smothness of $\partial \Omega$ is assumed, the boundary condition 
\[
\vc{v} \cdot \vc{n}|_{\partial \Omega} = 0
\]
is understood in the sense 
\[
\intO{ \Grad \phi \cdot \vc{v} } = 
- \intO{ \phi \Div \vc{v} }  
\]	
for any $\varphi \in C^1(\Ov{\Omega})$.	Similarly, the spaces $C^k(\Ov{\Omega})$ are defined as restrictions of functions in $C^k(R^d)$ to 
$\Ov{\Omega}$.
\end{Remark}		

We are ready to formulate our first result. 

\begin{Theorem} [\bf Discrete solutions] \label{mT1}
Let $\Omega \subset R^d$, $d=2,3$ be a bounded domain. 

There exists a set $I_{\vr,s}$ of the initial data 
\[
(\vr_0, s_0) \in I_{\vr,s}, 
\]
dense in $L^q(\Omega; (0, \infty)) \times L^q(\Omega)$, $1 \leq q <  \infty$, 
enjoying the following properties:	

\begin{itemize} 
\item For any $(\vr_0, s_0) \in I_{\vr,s}$, there exists $\vm_0 \in L^\infty(\Omega; R^d)$ such that the Euler system admits infinitely many 
admissible weak solutions $(\vr^\lambda, \vm^\lambda, s^\lambda)_{\lambda \geq 0}$ emanating from the initial data $(\vr_0, \vm_0, s_0)$.
\item For each $\lambda \geq 0$, there exists a finite integer $K = K(\lambda)$ such that 
\begin{equation} \label{m8}
	(\vr^\lambda, \vm^\lambda, s^\lambda) (t,x) \in \left\{ (\vr_k, \vm_k,  s_k) \Big| 
	k = 1, \dots, K \right\} \ \mbox{for a.a.}\ (t,x) \in (0,T) \times \Omega.
\end{equation}
for some constant vectors $(\vr_k, \vm_k, s_k) \in (0, \infty) \times R^d \times R$.	

\item 
\begin{align} 
	\intO{ \vr^{\lambda_1} s^{\lambda_1} (t, \cdot) } &\leq 
	\intO{ \vr^{\lambda_2} s^{\lambda_2} (t, \cdot) } \ \mbox{for a.a.}\ t \in (0,T),\br 
	\intO{ \vr^{\lambda_1} s^{\lambda_1}  } &\ne \intO{ \vr_{\lambda_2} s_{\lambda_2}  } \ \mbox{in}\ (0,T).
	\label{m9}
\end{align}
whenever $\lambda_1 < \lambda_2$.

\end{itemize}	
	
\end{Theorem}

\begin{Remark} \label{mR1}
	
It will be clear from the proof of Theorem \ref{mT1} that the solutions satisfy also the renormalized version of the entropy inequality, namely 	
\[
	\int_0^T \intO{ 
		\Big[ \vr Z(s) \partial_t \varphi + Z(s) \vm \cdot \Grad \varphi \Big] } \dt 
	\leq - \intO{ \vr_0 Z(s_0) \varphi (0, \cdot) } 
\] 
for any $\varphi \in C^1_c([0,T) \times \Ov{\Omega})$, $\varphi \geq 0$, 
and any $Z \in C(R)$, $Z'(s) \geq 0$.

\end{Remark}

\begin{Remark} \label{mR5}

As we shall see below, given $(\vr_0, s_0) \in I_{\vr, s}$, there exist infinitely many initial momenta $\vm_0$ satisfying Theorem \ref{mT1} with an arbitrarily large norm. In particular, we have 
\[
| \vm_0 | = \sqrt{ 2 \vr_0 \Lambda - \vr_0 d p(\vr_0, s_0) },  
\]
for any positive $\Lambda$, for which 
\[
2 \vr_0 \Lambda - \vr_0 d p(\vr_0, s_0) > 0.
\]
	
\end{Remark}	

Theorem \ref{mT1} claims the existence of infinitely many weak solutions 
for a dense set of initial densities and entropies. Our next result shows that there are infinitely many solutions 
of the Euler system with a prescribed initial density $\vr_0$ and the kinetic energy $E_{0, \rm kin}$. 
In particular, the initial kinetic energy need not be large as in Feireisl et al. \cite{FeKlKrMa}. 

\begin{Theorem}[\bf Discrete solutions with prescribed initial kinetic energy] \label{mT3} 
	Let $\Omega \subset R^d$, $d=2,3$ be a bounded domain.
	
There exists a set $I_{\vr, E_{\rm kin}}$ dense in 
$L^q(\Omega; (0, \infty)^2)$, $1 \leq q < \infty$ such that the following holds:
\begin{itemize}	
\item For any $(\vr_0, E_{0, {\rm kin}}) \in I_{\vr, E_{\rm kin}}$ there exist 
$(s_0, \vm_0) \in L^\infty(\Omega; R^{d+1})$,
\[
\frac{1}{2} \frac{ |\vm_0|^2}{\vr_0} = E_{0,{\rm kin}}\ \mbox{a.a. in}\ \Omega, 
\]
such that the Euler system admits infinitely many
admissible weak solutions $(\vr^\lambda, \vm^\lambda, s^\lambda)_{\lambda \geq 0}$ emanating from the initial data $(\vr_0, \vm_0, s_0)$.
\item For each $\lambda \geq 0$, there exists a finite integer $K = K(\lambda)$ such that 
\[
	(\vr^\lambda, \vm^\lambda, s^\lambda) (t,x) \in \left\{ (\vr_k, \vm_k, s_k) \Big| 
	k = 1, \dots, K \right\} \ \mbox{for a.a.}\ (t,x) \in (0,T) \times \Omega.
\]
for some constant vectors $(\vr_k, \vm_k, s_k) \in (0, \infty) \times R^d \times R$.	

\item 
\begin{align} 
	\intO{ \vr^{\lambda_1} s^{\lambda_1} (t, \cdot) } &\leq 
	\intO{ \vr^{\lambda_2} s^{\lambda_2} (t, \cdot) } \ \mbox{for a.a.}\ t \in (0,T),\br 
	\intO{ \vr^{\lambda_1} s^{\lambda_1}  } &\ne \intO{ \vr_{\lambda_2} s_{\lambda_2}  } \ \mbox{in}\ (0,T).
	\nonumber
\end{align}
whenever $\lambda_1 < \lambda_2$.

\end{itemize}	

\end{Theorem}	

Under certain restrictions, the conclusion of Theorem \ref{mT3} remains 
valid when we prescribe the initial density $\vr_0$ and the total energy 
$E_0$. 

\begin{Corollary} [\bf Discrete solutions with prescribed initial energy profile] \label{mC1}
Suppose $\Omega \subset R^d$ is a bounded domain, where either $d=2$, or $d=3$ and $\gamma < \frac{5}{3}$. Let $\Lambda > 0$ be a positive constant. 

There exists a set $I_{\vr, E}$ dense in 
$L^q(\Omega; (0, \infty)) \times L^q (0, \infty; J)$, 
where
\[
J = \left(\Lambda; \frac{2}{d} \frac{1}{\gamma - 1} \Lambda \right)
\]	
such that for any  $(\vr_0, E_0) \in I_{\vr, E}$ there exists 
$s_0 \in L^\infty(\Omega)$, $\vm_0 \in L^\infty(\Omega; R^d)$, 
\[
E_0 = \frac{1}{2} \frac{|\vm_0|^2}{\vr_0} + \vr_0 e(\vr_0, s_0),
\]
and the conclusion of Theorem \ref{mT3} holds.

\end{Corollary}

\begin{Remark} \label{rC1}
Note that thanks to the hypotheses in Corollary \ref{mC1}, 
\[
\frac{2}{d} \frac{1}{\gamma - 1} > 1 
\]
so the the open interval $J$ is non--empty. Moreover, it is easy to check that 
for any profile 
\[
E \in L^\infty(\Omega),\ \inf_{\Omega} E > 0
\]
there exists $\Lambda > 0$ such that 
\[
E + \Lambda \in \left( \Lambda; \frac{2}{d} \frac{1}{\gamma - 1} \Lambda \right).
\]	
\end{Remark}

\noindent
Theorems \ref{mT1}, \ref{mT3} as well as Corollary \ref{mC1} will be proved in 
Section \ref{e}.

Our final result asserts the existence of solutions with a given terminal profile $\widetilde{s}$. For simplicity, we restrict ourselves to the 
case $\gamma = \frac{5}{3}$ (a monoatomic gas), and $d=3$. Moreover, we suppose 
\begin{align} 
\underline{s} &\leq \widetilde{s}(x) \leq \Ov{s} \ \mbox{for all}\ x \in \Ov{\Omega}, \br	
\quad &\left| \left\{ x \in \Ov{\Omega} \ \Big| \ \widetilde{s} 
\ \mbox{is not continuous at}\ x \right\} \right|_d  = 0, 	
\label{m10}
\end{align}		
where $| \cdot |_d$ denotes the $d-$dimensional Lebesgue measure. In other words, the profile $\widetilde{s}$ belongs to the class of Riemann integrable 
functions in $\Ov{\Omega}$.

\begin{Theorem} [\bf Discrete solutions with prescribed terminal entropy profile] \label{mT2}
	Let $\gamma = \frac{5}{3}$ \\ (a monoatomic gas).
	Let $\Omega \subset R^3$ be a bounded domain. 
	
	Then there is a set $I_{\vr,s}$ of the initial data 
	\[
	(\vr_0, s_0) \in I_{\vr,s}, 
	\]
	dense in $L^q(\Omega; (0, \infty)) \times L^q(\Omega; (-\infty, \Ov{s}))$, $1 \leq q <  \infty$, 
	such that for any $(\vr_0, s_0) \in I_{\vr,s}$, and any 
	\[
	\Lambda >  \vr_0 e(\vr_0, s_0) ,
	\] 
	there exists $\vc{m}_0 \in L^\infty(\Omega; R^d)$ such that the following holds:
	
	\begin{itemize} 
		\item For an arbitrary entropy profile $\widetilde{s}$,
		\[
		s_0(x) < \widetilde{s} (x) ,\ \vr_0 e(\vr_0(x), \widetilde{s}(x)) < \Lambda 
		\ \mbox{for all}\ x \in \Ov{\Omega},\ \widetilde{s} \ \mbox{belonging to the class \eqref{m10}}, 
		\]
		the Euler system admits a weak solution $(\vr, \vm, s)$ emanating from 
		the initial data $(\vr_0, s_0,\vm_0)$ and satisfying 
		\begin{equation} \label{m11}
		\lim_{t \to T-} s(t, \cdot) = \widetilde{s} \ \mbox{strongly in}\ 
		L^q(\Omega)\ \mbox{for any}\ 1 \leq q < \infty.
		\end{equation}
		
		\item For any $0 < \tau < T$, there is an integer $K(\tau)$ such that
		\begin{equation} \label{m12}
			(\vr, \vm, s) (t,x) \in \left\{ (\vr_k, \vm_k, s_k) \Big| 
			k = 1, \dots, K(\tau) \right\} \ \mbox{for a.a.}\ (t,x) \in (0,\tau) \times \Omega.
		\end{equation}
	If, in addition, the profile $\widetilde{s}$ attains only a finite number of constant values, then \eqref{m12} holds for $\tau = T$.

	\end{itemize}	
	
\end{Theorem}	

\begin{Remark} \label{mR3}
Theorem \ref{mT2} holds for any $0 < T \leq \infty$.	
\end{Remark}	

\begin{Remark} \label{mR4}

The initial momentum $\vm_0$ depends only on $(\vr_0, s_0)$ and the 
constant $\Lambda$ but \emph{not} on the specific shape of the final profile $\widetilde{s}$. Moreover, it will be clear from the proof of 
Theorems \ref{mT1}, \ref{mT2} that the set $I_{\vr,s}$ of the initial data 
as well as the the momentum profile $\vm_0$ can be chosen the same in 
Theorems \ref{mT1}, \ref{mT2}. In particular, the same initial data give rise to the family of weak solution $(\vr^\lambda, \vm^\lambda, s^\lambda)_{\lambda \geq 0}$ claimed in Theorem \ref{mT1} as well as the  weak solutions converging to a prescribed final entropy profile $\tvr$ as in Theorem \ref{mT2}. 

\end{Remark}

The forthcoming three sections are devoted to the proof of Theorems \ref{mT1}, \ref{mT3}, and \ref{mT2}.

\section{Reformulation}
\label{r}

Our first task is to reformulate the original problem in the Convex Integration framework. The procedure has already been used in \cite{FeKlKrMa}. 

\subsection{The case of prescribed initial density and entropy}

Since the solutions claimed in Theorem \ref{mT1} attain only a discrete set 
of states, we focus on the 
initial data $\vr_0$, $s_0$ enjoying the same property.
Given a bounded domain $\Omega \subset R^d$, we consider its decomposition
into a finite union of $N$ subdomains,  
\begin{equation} \label{r1}
\Omega_n \subset \Omega \ \mbox{for all}\ n = 1, \dots, N, \ 
\Omega_i \cap \Omega_j = \emptyset \ \mbox{for}\ i \ne j,\ 
|\Omega \setminus \cup_{n=1}^N \Omega_n|_d = 0,
\end{equation}	
together with piecewise constant initial data, 
\begin{align}
\vr_0 \in L^\infty(\Omega),\ \vr_0(x) &= \vr_{0,n} > 0 \ \mbox{a positive constant in}\ \Omega_n,\br 
s_0 \in L^\infty(\Omega),\ s_0(x) &= s_{0,n} \in R \ \mbox{in}\ \Omega_n,\ n = 1, \dots, N.
\label{r2}
\end{align}	
The decomposition is arbitrary and no regularity of $\partial \Omega_n$ is required.
Obviously, the set of initial data $(\vr_0, s_0)$ belonging to the 
class \eqref{r1}, \eqref{r2} is dense 
in
$L^q(\Omega; (0, \infty)) \times L^q(\Omega)$ for any finite $1 \leq q <  \infty$, modulo a suitable refinement of the sets $\Omega_n$.

Following the strategy of \cite{FeKlKrMa} we set 
\begin{align} 
	\vr(t,x) &= \vr_{0,n} \ \mbox{for all}\ t \in [0,T] \times \Omega_n,\ 
	\br
		s(t,x) &= s_{0,n} \ \mbox{for all}\ t \in [0,T] \times \Omega_n,\ 
	n = 1,\dots, N.
	\label{r3}
\end{align}	

To determine the momentum $\vm$, we solve the problem: 
\begin{align}
\int_0^T \int_{\Omega_n} \vm_n \cdot \Grad \varphi \dx \dt &= 0 \br 
\mbox{for all}\ \varphi &\in C^1([0,T] \times \Ov{\Omega}_n), \label{r4} \\	 
\int_0^T \int_{\Omega_n} \left[ \vm_n \cdot \partial_t \bfphi + 
\left( \frac{\vm_n \otimes \vm_n}{\vr_{0,n}} - \frac{1}{d} \frac{|\vm_n|^2}{\vr_{0,n}} \mathbb{I}    \right) : \Grad \bfphi \right] \dx \dt &= 0 \br 
\mbox{for all}\ \bfphi &\in C^1([0,T] \times \Ov{\Omega}_n; R^d),
\label{r5}	
\end{align}	
supplemented with the kinetic energy constraint 
\begin{equation} \label{r6}
\frac{1}{2} \frac{ |\vm_n|^2 }{\vr_{0,n}} = \Lambda - \frac{d}{2} p(\vr_{0,n}, 
	s_{0,n}) \ \mbox{a.a. in}\ (0,T) \times \Omega_n, 
\end{equation}		
where $\Lambda$ is a suitable constant independent of $n=1,\dots, N$ chosen in such a way that
\begin{equation} \label{r7} 
\Lambda - \frac{d}{2} p(\vr_{0,n}, 
s_{0,n}) > 0 \ \mbox{for all}\ n = 1,\dots, N.
\end{equation}
Note carefully that the test functions in \eqref{r4}, \eqref{r5} do not vanish on the boundary of $\Omega_n$. In particular, the momentum defined as 
\begin{equation} \label{r8}
\vm (t,x) = \vm_n(t,x) \ \mbox{for}\ t \in (0,T),\ x \in \Omega_n ,\ 
n = 1,\dots, N,  	
\end{equation}	
satisfies the integral identities \eqref{r4}, \eqref{r5} with $\Omega_n$ replaced by $\Omega$.

It follows from \eqref{r2}, \eqref{r6} that $\vm_n \in L^\infty((0,T) \times \Omega_n; R^d))$. Moreover, by virtue of \eqref{r5}, 
\begin{equation} \label{r9}
\vm_n \in C_{\rm weak}([0,T] ; L^q(\Omega_n; R^d)) \ \mbox{for any finite}\ 
q > 1,\ \vm_n(0, \cdot) = \vm_n(T, \cdot) = 0 
\ \mbox{for all}\ n =1,\dots, N.
\end{equation}

The seemingly overdetermined system \eqref{r4}--\eqref{r6} can be solved via the method of Convex Integration elaborated in Section \ref{d}, cf. 
also \cite{FeKlKrMa}.

\subsection{The case of prescribed initial density and kinetic energy}

In order to adjust our ansatz to the hypotheses of Theorem \ref{mT3}, we go back to \eqref{r6} fixing a piecewise constant strictly positive function 
\begin{equation} \label{rr10}
E_{0,{\rm kin},n} = \frac{1}{2} \frac{|\vm_{n}|^2}{\vr_{0,n}} = 
\Lambda - \frac{d}{2} p(\vr_{0,n}, s_{0,n}) \ \mbox{in}\ \Omega_n.
\end{equation}
Given $E_{0, {\rm kin}, n}$, we may consider $\Lambda > E_{0,{\rm kin},n}$
for all $n = 1, \dots, N$ and adjust the initial entropy $s_{0,n}$ so that 
\begin{equation} \label{rr11}
\Lambda - E_{0,{\rm kin},n} = \frac{d}{2} p(\vr_{0,n}, s_{0,n}) 
\ \mbox{for}\ n = 1, \dots, N.
\end{equation}	

Similarly to the preceding section, we fix 
\begin{align} 
	\vr(t,x) &= \vr_{0,n} \ \mbox{for all}\ t \in [0,T] \times \Omega_n,\ 
	\br
	s(t,x) &= s_{0,n} \ \mbox{for all}\ t \in [0,T] \times \Omega_n,\ 
	n = 1,\dots, N, 
	\label{rr12}
\end{align}	
where $s_{0,n}$ are now determined by the prescribed 
initial kinetic energy via \eqref{rr11}. Thus the situation is the same as in the preceding section requiring solutions of \eqref{r4}--\eqref{r6} for $(\vr_{0,n}, s_{0,n})$ satisfying \eqref{rr11}.

\section{Discrete solutions via Convex Integration}
\label{d}

As already mentioned above, the method of Convex Integration provides \emph{infinitely many} solutions 
to problem \eqref{r4}--\eqref{r6}, see Chiodaroli \cite{Chiod}. Here, we use a refined version of this result due to Luo, Xie and Xin \cite{LuXiXi}. 
First we claim there exist $\frac{d(d+3)}{d}$ vectors $\vm^\ell_n \in R^d$, 
$\ell = 1, \dots, \frac{d(d+3)}{d}$ such that 
\begin{align}
\frac{1}{2} \frac{|\vm^{\ell}_n|^2}{\vr_{0,n}} &= 
\Lambda - \frac{d}{2} p(\vr_{0,n}, s_{0,n}), \label{r10} \\	
(0,0) &\in {\rm int} \left( {\rm convex \ hull} \left[ \left( \vm^{\ell}_n;  
\frac{\vm^{\ell}_n \otimes \vm^{\ell}_n}{\vr_{0,n}} - \frac{1}{d} 
\frac{|\vm^{\ell}_n|}{\vr_{0,n}} \right),\ \ell = 1, \dots, \frac{d(d+3)}{d} \right] \right), \br 
&\mbox{in}\ R^d \times R^{d \times d}_{0, {\rm sym}},  
\label{r11}
\end{align}
see \cite[Lemma 3]{LuXiXi}. The specific form of the vectors $\vm^\ell_n$ may depend on 
the constants $\vr_{0,n}$, $s_{0,n}$ as well as on $\Lambda$.

Now, we report the following crucial result shown by Luo, Xie, and Xin \cite[Proposition 1]{LuXiXi}. 

\begin{Proposition}[\bf Convex integration] \label{rP1}
Given $n=1,\dots, N$, problem \eqref{r4}, \eqref{r5} admits infinitely many solutions $\vm_n$ satisfying 
\begin{equation} \label{r12}
\vm_n(t,x) \in \left\{ \vm^\ell_n \Big| \ \ell = 1, \dots, \frac{d(d+3)}{2} \right\} \ \mbox{for a.a.}\ (t,x) \in (0,T) \times \Omega_n,
\end{equation}
with the constant vectors $\vm^\ell_n$ specified in \eqref{r10}, \eqref{r11}.	 
\end{Proposition}	

\begin{Remark} \label{RLL1}
	
More precisely, \cite[Proposition 1]{LuXiXi} is applied for 
$\underline{\vm} = 0$, 	$\underline{\vU} = 0$, $\vr = \vr_{0,n}$, $
q = 2 \vr_{0,n} \Lambda - d \vr_{0,n} {p(\vr_{0,n}, s_{0,n})}$, $\vc{B} \equiv 0$ in the notation 
of \cite{LuXiXi}.
	
\end{Remark}	

It is easy to see that \eqref{r10}, \eqref{r12} imply \eqref{r6}. In particular, equation \eqref{r5} can be written in the form 
\begin{equation} \label{r13}
\int_0^T \int_{\Omega_n} \left[ \vm_n \cdot \partial_t \bfphi + 
\left( \frac{\vm_n \otimes \vm_n}{\vr_{0,n}} \right) : \Grad \bfphi +
p(\vr_{0,n}, s_{0,n}) \Div \bfphi - \frac{2}{d}\Lambda \Div \bfphi \right] \dx \dt = 0 
\end{equation} 
for all $\bfphi \in C^1([0,T] \times \Ov{\Omega}_n; R^d)$. Thus summing up 
over $n$ and using test functions with vanishing normal trace, we conclude 
\begin{equation} \label{r14}
	\int_0^T \int_{\Omega} \left[ \vm \cdot \partial_t \bfphi + 
	\left( \frac{\vm \otimes \vm}{\vr} \right) : \Grad \bfphi +
	p(\vr, s) \Div \bfphi \right] \dx \dt = 0 
\end{equation}
for all $\bfphi \in C^1([0,T] \times \Ov{\Omega}; R^d)$, $\bfphi \cdot \vc{n}|_{\partial \Omega} = 0$.

In view of \eqref{r4}, it is a routine matter to check that 
\begin{equation} \label{r15}
	\int_0^T \intO{ \Big[ \vr \partial_t \varphi + \vm \cdot \Grad \varphi \Big] } \dt = 
	- \intO{ \vr_0 \varphi(0,\cdot)} 
\end{equation}	
for any $\varphi \in C^1_c([0,T) \times \Ov{\Omega})$, and 
\begin{equation} \label{r16}
	\int_0^T \intO{ 
		\Big[ \vr s \partial_t \varphi + s \vm \cdot \Grad \varphi \Big] } \dt 
	= - \intO{ \vr_0 s_0 \varphi (0, \cdot) } 
\end{equation} 
for any $\varphi \in C^1_c([0,T) \times \Ov{\Omega})$. Thus the entropy balance is satisfied as equality. 

Finally, seeing that the total energy 
\[
E = \frac{1}{2} \frac{|\vm|^2}{\vr} + \vr e (\vr, s) 
\]
is piecewise constant, specifically, 
\[
E = \Lambda + \vr_{0,n} e(\vr_{0,n}, s_{0,n}) - \frac{d}{2} p(\vr_{0,n}, s_{0,n}) \ \mbox{a.a. in}\ (0,T) \times \Omega_n, \ n =1, \dots, N,
\]
we may repeat the above arguments to deduce the energy balance 
\begin{equation} \label{r17}
	\int_0^T \intO{ \left[ E(\vr, \vm,s ) \partial_t \varphi + 
		\Big( E(\vr, \vm, s) + p(\vr,s) \Big) \frac{\vm}{\vr} \cdot \Grad \varphi \right]	} \dt = 0 
\end{equation}				
for any $\varphi \in C^1_c((0,T) \times \Ov{\Omega})$. Finally, seeing that 
the total energy is constant, we may extend \eqref{r17} to 
\begin{equation} \label{r17a}
	\int_0^T \intO{ \left[ E(\vr, \vm,s ) \partial_t \varphi + 
		\Big( E(\vr, \vm, s) + p(\vr,s) \Big) \frac{\vm}{\vr} \cdot \Grad \varphi \right]	} \dt = 0 
\end{equation}
for any $\varphi \in C^1([0,T] \times \Ov{\Omega})$, $\varphi(0, \cdot) = 
\varphi (T, \cdot)$.

At this stage, one is tempted to conclude that $(\vr, s, \vm)$ is an admissible 
weak solution of the Euler system. However, this is not true since, by virtue of \eqref{r14}, 
\[
\vm \in C_{\rm weak} ([0,T]; L^q(\Omega; R^d)) 
\ \mbox{for any}\ 1 \leq q < \infty,\ \vm(0, \cdot) = \vm(T, \cdot) = 0,
\]
while the initial energy 
\[
E(0, \cdot) = \Lambda + \vr_0 e(\vr_0, s_0) - \frac{d}{2}p (\vr_0, s_0) > 
\vr_0 e(\vr_0, s_0). 
\]
In other words, the (total) energy experiences a positive jump at the 
initial time $t = 0$ - a phenomenon pertinent to solutions obtained via 
Convex Integration. 

Fortunately, by virtue of \eqref{r6}, we have  
\[
E(\vr, s, \vm)(t, \cdot) = \frac{1}{2} \frac{|\vm (t, \cdot) |^2}{\vr} 
+ \vr e(\vr, s) (t, \cdot) 
\]
for a.a. $t \in (0,T)$.
Thus performing a simple time shift, we have proved the following result. 

\begin{Proposition}[\bf Discrete solutions] \label{rP2}
Let the domain decomposition \eqref{r1}, the discrete initial data $(\vr_{0,n})_{n=1}^N$, $(s_{0,n})_{n=1}^N$, and the constant $\Lambda$ satisfying \eqref{r7} be given. 

Then for a.a. $\tau \in (0, T)$ there exists an initial momentum $\vm_0 \in L^\infty$ such that 
the Euler system \eqref{i1}--\eqref{i8} admits an admissible weak solution $(\vr, s, \vm)$ in $(0,\tau) \times \Omega$ satisfying:
\begin{itemize}
\item	
\begin{align}
\vr &\in L^\infty((0,\tau) \times \Omega),\ 	
\vr(t,x) = \vr_{0,n} \ \mbox{for} \ t \in [0,\tau],\ x \in \Omega_n, \br s &\in L^\infty((0,\tau) \times \Omega),\ s(t,x)= s_{0,n} \ \mbox{for}\ t \in (0, \tau), \ x \in \Omega_n,\ n=1,\dots, N;	 
\label{r18} 
\end{align}
\item
\begin{align}
 \vm &\in L^\infty((0, \tau) \times \Omega; R^d) \cap C_{\rm weak}([0,\tau]; L^q(\Omega; R^d)),\ 
1 \leq q < \infty, \br 
\vm(0, \cdot) &= \vm_0,\ \vm(\tau, \cdot) = 0,
\label{r19} 
\end{align}
\item
\begin{align}
\vm &= \vm_n(t,x) \in \left\{ \vm^{\ell}_n\Big|\ \ell = 1, \dots, \frac{d(d+3)}{2} \right\} \ \mbox{for a.a.}\ (t,x) \in (0,\tau) \times \Omega_n,
\label{r20}
\end{align}
\item
\begin{align}
E &\in L^\infty((0,\tau) \times \Omega),\br E(t,x) &= E_n = \frac{1}{2} \frac{|\vm_n|^2}{\vr_{0,n}} + \vr_{0,n} e(\vr_{0,n} s_{0,n}) = 
\Lambda + 	\vr_{0,n} e(\vr_{0,n} s_{0,n}) - \frac{d}{2} p(\vr_{0,n} s_{0,n}), \br
&\mbox{for a.a.}\ (t,x) \in (0,\tau) \times \Omega_n.
\label{r21}
\end{align}		
\end{itemize}
\end{Proposition}	

Note carefully that \eqref{r21} in the context of Theorem \ref{mT3} reads 
\[
\frac{1}{2} \frac{|\vm_n|^2}{\vr_{0,n}} (t,x) = \Lambda - \frac{d}{2} p(\vr_{0,n}, s_{0,n}) = E_{0, {\rm kin}, n}\ \mbox{for a.a.}\ (t,x) \in (0,\tau) 
\times \Omega_n.
\]

All components of the weak solutions obtained in Proposition \ref{rP2} are constant in the space--time cylinders $(0,\tau) \times \Omega_n$, $n=1,\dots, N$, with the exception of the momentum $\vm_n$. However, the modulus $|\vm_n|^2$ proportional to the kinetic energy enjoys the same property as well. In addition, in accordance 
with \eqref{r20}, the momentum $\vm$ attains in $\Omega_n$ only a finite number of constant values all of them on the sphere of the diameter
\[
|\vm_n| = \sqrt{ 2 \vr_{0,n} \Lambda - d \vr_{0,n}p(\vr_{0,n}, s_{0,n})}.
\]

We point out that Proposition \ref{rP2} guarantees the existence of only 
\emph{one} discrete solution. We believe that elaborating more carefully the technique developed in \cite{LuXiXi}, we could obtain the existence of infinitely many solutions already at this level. Note, however, that all these solutions would \emph{conserve} entropy in contrast with the claim of Theorems \ref{mT1}, \ref{mT3}, \ref{mT2}.

\section{Entropy production rate}
\label{e}

We are ready to complete the proofs of our main results claimed in Theorems \ref{mT1}, \ref{mT3}, and \ref{mT2}. 

\subsection{Proof of Theorems \ref{mT1}, \ref{mT3} - solutions with a prescribed entropy production rate and initial kinetic energy}

 The idea is quite simple; we divide the time interval $(0,T)$ into a finite number of segments and augment the entropy on each of them. Note that such a process is compatible with the entropy \emph{inequality} \eqref{m7}.

Let 
\[
0 = \tau_0 < \tau_1 < \dots < \tau_M = T,\ T \ \mbox{finite},\ 
\tau_M < \infty \ \mbox{if}\ T = \infty. 
\]
be a division of the time interval $[0,T)$, where 
$\tau_1$ belongs to the subset of full measure in $(0,T)$, for which Proposition 
\ref{rP2} guarantees the existence of a weak solution in the class 
\eqref{r18}--\eqref{r21}.
Let 
the initial data $\vr_0$, $s_0$, and $\vm_0$ as well as the constant $\Lambda$ 
be the same as in Proposition \ref{rP2}. Applying Proposition \ref{rP2}, we obtain a weak solution of the Euler system in the space time domain $(0,\tau_1) \times \Omega$ satisfying \eqref{r18}--\eqref{r21}.

Now, we repeat the whole construction elaborated in Section \ref{d} on the 
time interval $(\tau_1, \tau_2)$, with the initial data 
\begin{align} 
\vr(\tau_1, \cdot) &= \vr_0,\  \vr(\tau_1, \cdot) = \vr_{0,n}\ \mbox{in}\ \Omega_n, \ n=1,\dots, N, \br
\vm(\tau_1, \cdot) &= 0,
\label{e1}
\end{align}
and 
\begin{equation} \label{e2}
s(\tau_1, \cdot) = \widehat{s},\ \ \mbox{for some}\  
\widehat{s} = \widehat{s}_{0,n} \ \mbox{in}\ \Omega_n,\ 
\widehat{s}_{0,n} > s_{0,n} \ \mbox{in}\ \Omega_n.	
\end{equation}
As pointed out above, the entropy profile $\widehat{s}$ is compatible with the 
entropy inequality.

Finally, we fix a constant $\widehat{\Lambda}$ satisfying 
\begin{equation} \label{e3} 
	\widehat{\Lambda} - \frac{d}{2} p(\vr_{0,n}, 
	\widehat{s}_{0,n}) > 0 \ \mbox{for all}\ n = 1,\dots, N.
\end{equation}

The arguments specified in Section \ref{d} yield the existence of a weak solution on the 
time interval $(\tau_1, \tau_2)$ enjoying the properties \eqref{r18}--\eqref{r21}, and with the initial momentum 
\[
\vm(\tau_1, \cdot) = 0.
\]
In particular, by virtue of \eqref{r19}, the two solutions can be concatenated 
at the time $\tau_1$
yielding a new solution in $(0,\tau_2) \times \Omega$ 
at level of the conservation mass, momentum, and entropy production as 
$\vm$ is weakly continuous at $t = \tau_1$. The new solution satisfies 
the entropy inequality \eqref{m7}, with the entropy production rate 
\[
\sigma = \delta_{t = \tau_1} \vr (\widehat{s} - s) > 0,
\]
specifically 
\[
\partial_t (\vr s) + \Div (s \vm) = \sigma \geq 0 
\]
in the sense of distributions in $(0, \tau_2) \times \Omega$.

Finally, we have to make sure that the energy conservation \eqref{m6} holds in  the whole interval $(0, \tau_2)$, meaning,  the total energy must 
equal the same constant on each set $\Omega_n$, $n=1, \dots, N$. In accordance with \eqref{r21}, this amounts to showing 
\begin{equation} \label{e4}
\Lambda + \vr_{0,n} e(\vr_{0,n}, s_{0,n}) - \frac{d}{2} 
p(\vr_{0,n}, s_{0,n}) = \widehat{\Lambda} + \vr_{0,n} e(\vr_{0,n}, \widehat{s}_{0,n}) - \frac{d}{2} 
p(\vr_{0,n}, \widehat{s}_{0,n}) \ \mbox{for all}\ n = 1, \dots, N.
\end{equation}
Moreover, by virtue of \eqref{m1}, this is equivalent to 
\begin{equation} \label{e5}
	\Lambda - \widehat{\Lambda} = 
	\left(1 - \frac{d}{2}(\gamma - 1) \right) \vr_{0,n} \Big( 
	e(\vr_{0,n}, \widehat{s}_{0,n}) - e(\vr_{0,n}, {s}_{0,n}) \Big)
	\ \mbox{for all}\ n = 1, \dots, N,
\end{equation}	
where 	$\left(1 - \frac{d}{2}(\gamma - 1) \right) \geq 0$ as long as $1 < \gamma \leq \frac{2}{d} + 1$.
In view of Gibbs' relation \eqref{i4}, the internal energy $e$ is a strictly 
increasing function of the entropy $s$, and 
there are two alternatives: 
\begin{enumerate}
	\item $d=3$, $\gamma = \frac{5}{3}$. We may consider $\widehat{\Lambda}= 
	\Lambda$, and the entropy in the form
	\[ 
	\widehat{s}_{0,n} = {s}_{0,n} + \lambda,\ \lambda \geq 0 
	\]
as long as the parameter $\lambda$ satisfies	
	\[
	\Lambda > \frac{3}{2} p(\vr_{0,n}, s_{0,n} + \lambda). 
	\]
	\item 
	$\widehat{\Lambda} = \Lambda - \lambda$, $\lambda > 0$, where $\lambda$ is a small positive parameter. The entropy $\widehat{s}$ is then chosen as
	\[
	\widehat{s}_{0,n} = \widehat{s}_{0,n}(\lambda), 
	\]
	solving the equation
	\[
	\lambda = 
	\left(1 - \frac{d}{2}(\gamma - 1) \right) \vr_{0,n} \Big( 
	e(\vr_{0,n}, \widehat{s}_{0,n} (\lambda) ) - e(\vr_{0,n}, {s}_{0,n}) \Big)
	\ \mbox{for all}\ n = 1, \dots, N,
	\]
under the constraint 
\[
\widehat{\Lambda} > \frac{3}{2} p(\vr_{0,n}, \widehat{s}_{0,n}). 
\]

\end{enumerate}	

In both cases, the procedure yields a family of solutions $(\vr^\lambda, \vm^\lambda, s^\lambda)_{\lambda \geq 0}$ claimed in Theorems \ref{mT1}, \ref{mT3}. The process can be continued up to $\tau_M = T$. 
We have proved Theorems \ref{mT1}, \ref{mT3}.

To show Corollary \ref{mC1} it is enough to observe that 
\begin{equation} \label{for}
E = \frac{|\vm|^2}{\vr} + \vr e(\vr,s) = 
\Lambda - \frac{d}{2} p(\vr, s) + \vr e(\vr,s) = 
\Lambda + \left( \frac{1}{\gamma - 1} - \frac{d}{2} \right) p(\vr,s).
\end{equation}
Since 
\[
 \left( \frac{1}{\gamma - 1} - \frac{d}{2} \right) > 0,
\]
the discrete values of the total energy $E_{0,n}$ in $\Omega_n$ can be prescribed in terms $(\vr_{0,n}, s_{0,n})$ as long as $E_{0,n} > \Lambda$, 
and the kinetic energy remains positive, meaning 
\[
\frac{d}{2} p(\vr_{0,n}, s_{0,n}) < \Lambda \ \Rightarrow \ 
E_{0,n} < \frac{2}{d} \frac{1}{\gamma - 1} \Lambda.
\]

Clearly, the above construction can be repeated for a different choice of the 
times $(\tau_i)_{i=1}^I$ to obtain a large variety of solutions. For instance, 
one could asymptotically approach a given total entropy profile $\Phi$,
\[
t \in [0,T] \mapsto \intO{ \vr s (t,x) } = \Phi(t)
\] 
In particular, the time $\tau_1$ can be chosen arbitrarily close to $t = 0$. 
This observation yields the following corollary. 

\begin{Corollary} \label{eC1}
	
Given $\ep > 0$, the solution family $(\vr^\lambda, \vm^\lambda, s^\lambda)_{\lambda \geq 0}$ claimed in Theorems \ref{mT1}, \ref{mT3} can be constructed in such a way that 
\[
(\vr^{\lambda_1}, \vm^{\lambda_1}, s^{\lambda_1}) \ne 
(\vr^{\lambda_2}, \vm^{\lambda_2}, s^{\lambda_2}) \ \mbox{in} \ 
(0, \ep) \times \Omega \ \mbox{for}\ \lambda_1 \ne \lambda_2.
\]	
In particular, the initial data $(\vr_0, s_0, \vm_0)$ identified in Theorems 
\ref{mT1}, \ref{mT3} are \emph{wild} in the sense of \cite{ChFe2023}. 
\end{Corollary}

Finally, we note that, by virtue of Proposition \ref{rP1}, the problem on each interval $(\tau_m, \tau_{m+1})$, $1 \leq m < M$,  
admits infinitely many solutions. 

\subsection{Proof of Theorem \ref{mT2} - solutions with a prescribed terminal entropy profile}		

Similarly to the preceding part, we consider an (infinite) sequence of times
\[
0 = \tau_0 < \tau_1 < \dots \tau_m < T,\ t_m \to T \ \mbox{as}\ m \to \infty,
\]
and concatenate solutions at the points $\tau_m$, m=1,2,\dots

Similarly to \eqref{e1}, we fix the initial data 
\[
\vr(\tau_m, \cdot) = \vr_0,\ \vr(\tau_m, \cdot) = \vr_{0,n}\ \mbox{in}\ \Omega_n,\ n =1, \dots, N,\ m = 0, 1, \dots, 
\]
and 
\[
\vm(\tau_n, \cdot) = 0,\ m =1,2,\dots. 
\]
In addition, keeping in mind $d = 3$, $\gamma = \frac{5}{3}$, 
in particular, 
\[
p = \frac{2}{3} \vr e,
\]
we choose $\Lambda > 0$ so that 
\[
\Lambda - \frac{3}{2} p(\vr_{0,n}, s_{0,n}) > 0,\ n = 1, \dots, N  
\]
 Recalling formula \eqref{e4} and the hypothesis $d = 3$, $\gamma = \frac{5}{3}$, the constant $\Lambda$, and consequently the total energy will be the same at any time interval $\tau_m, \tau_{m+1}$ as long as the initial entropies are chosen so the $s(\tau_m, \cdot) \leq \Ov{s}$. 

Finally, we define the entropy ``initial data'' at $\tau_1, \tau_2, \dots$ by recursion. 
In accordance with the hypotheses of Theorem \ref{mT2}, we start with 
\[
s_0 (x) = s_{0,n} \ \mbox{for}\ x \in \Omega_n \equiv \Omega_{0,n},\ 
n = 1,\dots, N \equiv N(0),
\]
\[
s_0(x) < \widetilde{s} (x),\ \frac{3}{2} p(\vr_{0}(x), \widetilde{s}(x)) < \Lambda.
\]

Next, for $m=1, 2, \dots$, we define recursively a refinement of 
$\Omega$-decomposition:
\begin{itemize}
\item 
\begin{align}
\Omega_{m,n} &\subset \Omega,\ m = 0,1,\dots ,\ n = 1, \dots, N(m),\ 
\Omega_{m,i} \cap \Omega_{m,j} = \emptyset \ \mbox{for}\ i \ne j, \br  
| \Omega \setminus \cup_{n=1}^{N(m)} \Omega_{m,n}|_d &= 0 \ \mbox{for}\ m =1,2, \dots; 
\nonumber
\end{align}
\item
For any $\Omega_{m,n}$, there is $\Omega_{m-1,j}$ such that 
\begin{equation} \label{e6}
\Omega_{m,n} \subset \Omega_{m-1,j}
\end{equation}
for any $m = 1,2,\dots$;
\item
\begin{equation} \label{e7}
\sup_{n = 1,\dots, N(m)} {\rm diam}[\Omega_{m,n}] \to 0 \ \mbox{as}\ m \to \infty.	
\end{equation}
\end{itemize}

Finally, we fix the initial data for the entropy, 
\[
s(\tau_m, \cdot) = \inf_{y \in \Omega_{m,n}} \widetilde{s}(y)
\ \mbox{for}\ x \in \Omega_{m,n},\ n=1, \dots, N(m).
\]

Now, thanks to \eqref{e6}, we have 
\[
s(\tau_{m-1}, \cdot) \leq s(\tau_{m}, \cdot) \ \mbox{for all}\ 
m=1,2 ,\dots.
\]
In particular, our choice of piecewise constant initial entropies is compatible 
with the entropy inequality \eqref{m7}. We also note that $\vr(\tau_m, x) = 
\vr_{0,n}$, whenever $x \in \Omega_n$, in particular the density remains continuous.

Finally, we have 
\[
\intO{ s(t, \cdot)} = \sum_{n=1}^{N(m)} (\inf_{\Omega_{n,m}} \widetilde{s}) 
|\Omega_{m,n}| \ \mbox{for}\ t \in [\tau_m, \tau_{m+1}), 
\]
where, by virtue  of \eqref{e7}, 
\[
\sum_{n=1}^{N(m)} (\inf_{\Omega_{n,m}} \widetilde{s}) 
|\Omega_{m,n}| \to \intO{ \widetilde{s} } 
\]
as the profile $\widetilde{s}$ is Riemann integrable. Seeing that the sequence 
$(s(\tau_m, \cdot))_{m=1}^\infty$ is non--decreasing, we conclude 
\[
\lim_{t \to T-} \| s(t, \cdot) - \widetilde{s} \|_{L^1(\Omega)} = 0.
\]

We have proved Theorem \ref{mT2}.

\section{Discrete weak solutions, admissibility, and the long time behaviour}
\label{C}

We finish the paper by discussing the implications of the results in the 
broader context of the theory of general dissipative measure--valued solutions 
of the complete Euler system, see e.g. \cite[Part II, Chapter 5]{FeLMMiSh}, 
in particular their admissibility and the long time behaviour. 
As confirmed by the results of the present paper, the Euler system in 
the physically relevant multidimensional case is ill posed even in the class of 
entropy admissible weak solutions.

\subsection{General equations of state}

The polytropic constraint \eqref{m1} may be restrictive for certain models arising in astrophysics, where the pressure includes the effect of thermal radiation: 
\[
p(\vr, \vt) = \vr \vt + a \vt^4, 
\]
see e.g. Battaner \cite{BATT}. Revisiting the proof of Theorem \ref{mT1}, \ref{mT3}, specifically relation \eqref{e5}, we easily observe that it can be extended to the case 
\[
p(\vr,s) = \sum_{i=1}^I p_i(\vr,s),\ e(\vr,s) = \sum_{i=1}^I e_i(\vr,s),\ 
p_i(\vr,s) = (\gamma_i - 1) \vr e(\vr,s),
  \]
\[
1 < \gamma_i \leq \frac{2}{d} + 1,\ i = 1,\dots, I.
\] 

More general equations of state including a non--monotone dependence of $p$ 
on the density $\vr$ can be handled case by case as long as Gibbs' equation \eqref{i4} is satisfied.

\subsection{Measure valued solutions}

Any weak solution $(\vr, \vm, s)$ can be associated with a mono--atomic measure valued solution 
\[
(t,x) \in (0,T) \times \Omega \to \delta_{(\vr, \vm, s)(t,x)} \in 
\mathfrak{P}(R^{d+2}),
\]
where the symbol $\mathfrak{P}$ stands for the family of Borel probability measures. Moreover, it is easy to check that any \emph{convex combination} 
of measure--valued solutions is a measure--valued solution (cf. \cite[Part II]
{FeLMMiSh}). Consequently, the quantity 
\[
(t,x) \in (0,T) \times \Omega \mapsto \sum_{i = 1}^I \mu_i \delta_{(\vr^{\lambda_i}, \vm^{\lambda_i}, s^{\lambda_i})(t,x)}, \ 
\mu_i \geq 0,\ \sum_{i=1}^I \mu_i = 1
\] 
where $(\vr^{\lambda_i}, \vm^{\lambda_i}, s^{\lambda_i})$ are the solutions 
claimed in Theorems \ref{mT1}, \ref{mT3}, provides an example of a (truly) measure--valued solution attaining only a finite number of constant values.

\subsection{Admissibility}

Several attempts have been made to restore well--posedness of the Euler system by imposing suitable \emph{admissibility criteria} that go beyond the mere satisfaction of the entropy inequality \eqref{i6}. 

As a matter of fact, \eqref{i6} can be interpreted as 
\begin{equation} \label{C1}
\partial_t (\vr s) + \Div ( s \vm) = \sigma \geq 0, 	
\end{equation}
where $\sigma$ can be interpreted as a non--negative measure representing the 
\emph{entropy production rate}.
There are two selection criteria based on the specific form of $\sigma$: 

In his pioneering paper \cite{DiP2}, DiPerna introduced the following partial 
ordering in the class of solutions of the Euler system,  
\begin{align}
	(\vr^1, \vm^1, s^1) &\prec_{DiP} (\vr^2, \vm^2, s^2) \br
	\qquad &\Leftrightarrow_{\rm def}  \qquad \br 
	\intO{ \vr^1 s^1 (t+, \cdot) } &\leq  	\intO{ \vr^2 s^2 (t+, \cdot) } \ \mbox{for all}\ t \in (0,T). 
	\nonumber
\end{align}	
DiPerna conjectured that the \emph{measure--valued} solutions of the Euler system that are maximal with respect to $\prec_{DiP}$ must be weak solutions. 
In view of Theorem \ref{mT1}, in particular \eqref{m9}, we conclude that 
the solutions claimed in Theorem \ref{mT1} \emph{are not} maximal with respect 
to $\prec_{DiP}$. 

Dafermos \cite{Dafer} proposed a refined selection criterion based on partial ordering 
\begin{align}
	(\vr_1, \vm_1, s_1) &\prec_{Daf}
	(\vr_2, \vm_2, s_2)\br \ &\Leftrightarrow_{\rm{def}} \br
	\mbox{there exists}\ \tau \geq 0 \ \mbox{such that}\ (\vr_1, \vm_1, s_1)(t, \cdot) &=
	(\vr_2, \vm_2, s_2)(t, \cdot) \ \mbox{for all}\
	t \leq \tau, \br
	\frac{\D^+}{\dt} \intO{ \vr_2 s_2 (\tau, \cdot) } &>
	\frac{\D^+}{\dt}  \intO{ \vr_1 s_1  (\tau, \cdot) }.
	\nonumber
\end{align}
As we have shown in \cite{FeiLukYu}, the measure valued solutions that are maximal with respect to $\prec_{Daf}$ must be weak solutions. Unfortunately, the weak solutions claimed in Theorems \ref{mT1}, \ref{mT3} \emph{are not} 
maximal with respect to $\prec_{Daf}$. This can be seen in the proof 
based on the choice of the time interval decomposition 
$\tau_0 = 0 < \tau_1 < \dots \tau_M = T$. Indeed we can always construct a solution based on a new decomposition $\tau'_0 = 0 < \tau'_1 < \dots $, where 
$\tau'_1 < \tau_1$.

Like most of their counterparts obtained via Convex Integration, the 
``discrete'' solutions discussed in the present paper are eliminated by selection criteria based on the maximality of the entropy production rate. Very roughly indeed, one may say that the same method can be used to produce a new solution that is ``greater'' than the original one. This fact can be interpreted in a positive way that these solutions might be \emph{unphysical}. 
Note that the selection criteria based on maximality of the entropy production 
can identify a unique semigroup solution of the Euler system, see 
\cite{BreFeiHof19C}, \cite{FeiLukYu}.

\subsection{Long--time behaviour}

The weak solutions obtained in the present paper feature several 
``pathological'' properties in the asymptotic limit $t \to \infty$. To begin with,  the initial density profile $\vr_0$ is conserved. Solutions obtained in 
Theorem \ref{mT1} conserve also the entropy profile $s_0$ or its slight modification after the time $\tau_1$. Solutions claimed in Theorem \ref{mT3}
conserve even the initial distribution of the kinetic energy 
$
\frac{1}{2} \frac{ |\vm|^2 }{\vr}
$
or modification after $\tau_1$.
These properties evidently suggest that the solutions are not physically realistic. Still, they are entropy admissible solutions of the Euler system of gas dynamics.

\subsubsection{Time--periodic solutions}

It is worth noting that the weak solutions obtained in Section \ref{d} via 
Proposition \ref{rP1} are, in fact, \emph{time periodic} entropy admissible weak solutions of the Euler system, cf. the energy balance \eqref{r17a}.
Accordingly, we may deduce the following form of the result stated in Theorem 
\ref{mT1}.

\begin{Theorem}[\bf Time periodic solutions] \label{CT1}
	
Let $\Omega \subset R^d$, $d=2,3$ be a bounded domain. 	

There exists a set $I_{\vr,s}$ of the initial data 
\[
(\vr_0, s_0) \in I_{\vr,s}, 
\]
dense in $L^q(\Omega; (0, \infty)) \times L^q(\Omega)$, $1 \leq q <  \infty$, 
enjoying the following properties:	
\begin{itemize}
\item For any $(\vr_0, s_0) \in I_{\vr,s}$, there exists $\vm_0 \in L^\infty(\Omega; R^d)$ such that the Euler system admits infinitely many 
admissible weak solutions $(\vr^T, \vm^T, s^T)_{T > 0}$ emanating from the initial data $(\vr_0, \vm_0, s_0)$.
\item For each $T > 0$, there exists a finite integer $K = K(T)$ such that 
\[
	(\vr^T, \vm^T, s^T) (t,x) \in \left\{ (\vr_k, \vm_k, s_k) \Big| 
	k = 1, \dots, K \right\} \ \mbox{for a.a.}\ (t,x) \in (0,T) \times \Omega.
\]
for some constant vectors $(\vr_k, \vm_k, s_k) \in (0, \infty) \times R^d \times R$.
\item The solutions are time-periodic with the period $T>0$ starting from a positive time, say, $\tau_1 > 0$:
\begin{equation} \label{CC1}
(\vr^T, \vm^T, s^T)(t + T, \cdot) = 
(\vr^T, \vm^T, s^T)(t, \cdot) \ \mbox{for any}\ t > \tau_1.
\end{equation}
\end{itemize}
\end{Theorem}	

\begin{Remark} \label{CR1}
A similar result can be shown under the hypotheses of Theorem \ref{mT3}. We leave the details to the interested reader.	
\end{Remark}

It worth noting that the periodicity stated in \eqref{CC1} indeed holds for 
\emph{any} $t > \tau_1$ as $\vr^T$, $s^T$ are strongly continuous 
in $[\tau_1, \infty)$, while $\vm^T$ is weakly continuous in $[\tau_1, \infty)$ 
ranging in $L^q(\Omega)$, $1 \leq q < \infty$.

\subsubsection{Solutions converging to equilibria}

Revisiting the proof of Theorem \ref{mT3}, we can show the existence of global in time weak solutions converging to equilibrium $\vm = 0$. For the sake of simplicity, we follow the hypotheses of Theorem \ref{mT2} focusing on the case 
$d = 3$, $\gamma = \frac{5}{3}$. Evoking the situation 
handled in the proof of Theorem \eqref{mT1} we choose a sequence of time 
$\tau_0 < \tau_1 < \tau_2 < \dots < \tau_m \to \infty$. 
At each point $\tau_m$ we choose the initial data for the entropy to be 
\[
s(\tau_m , x) = s_{m,n} \ \mbox{for}\ x \in \Omega_n, \ 
s_{m,n} \geq s_{m-1,n},  
\]
\[
\Lambda - \frac{3}{2} p(\vr_{n}, s_{m,n}) \searrow 0 \ \mbox \ \mbox{as}\ 
m \to \infty \ \mbox{uniformly for}\ n = 1, \dots, N.
\]
Note that the quantity 
\[
\Lambda - \frac{3}{2} p(\vr_{n}, s_{m,n}) = \frac{1}{2} \frac{|\vm(x)|^2}{\vr_n(x)},\ x \in \Omega_n.  
\]
represent the kinetic energy of the solution in the time interval 
$(t_m, t_{m+1})$. As the choice of 
\[ 
\Lambda - \frac{3}{2} p(\vr_{n}, s_{m,n}) 
\]
and the sequence $\tau_m \to \infty$ was arbitrary, we have shown the following result. 

\begin{Theorem}[\bf Convergence to equilibria] \label{CT2}
	
	Let $\Omega \subset R^3$, be a bounded domain, and let $\gamma = \frac{5}{3}$. 	
	
	There exists a set $I_{\vr,s}$ of the initial data 
	\[
	(\vr_0, s_0) \in I_{\vr,s}, 
	\]
	dense in $L^q(\Omega; (0, \infty)) \times L^q(\Omega)$, $1 \leq q <  \infty$, 
	enjoying the following properties:	
	\begin{itemize}
		\item For any $(\vr_0, s_0) \in I_{\vr,s}$, there exists $\vm_0 \in L^\infty(\Omega; R^d)$ such that the Euler system admits infinitely many 
		admissible weak solutions $(\vr^\lambda, \vm^\lambda, s^\lambda)_{\lambda > 0}$ in $(0, \infty) \times \Omega$ emanating from the initial data $(\vr_0, \vm_0, s_0)$.
		\item For each $T > 0$, there exists a finite integer $K = K(T, \lambda)$ such that 
		\[
		(\vr^\lambda, \vm^\lambda, s^\lambda) (t,x) \in \left\{ (\vr_k, \vm_k, s_k) \Big| 
		k = 1, \dots, K \right\} \ \mbox{for a.a.}\ (t,x) \in (0,T) \times \Omega.
		\]
		for some constant vectors $(\vr_k, \vm_k, s_k) \in (0, \infty) \times R^d \times R$.
		\item For any profile 
		\[
		\chi \in C^1[0, \infty),\ \chi(t) > 0,\ \chi'(t) \leq 0,\ \chi(t) \to 
		0 \ \mbox{as} \ t \to \infty,
		\] 
		there exist $\underline{\lambda} = \underline{\lambda}(\chi)$, $\Ov{\lambda} = \Ov{\lambda}(\chi)$, and a time $\tau = \tau(\chi) > 0$ such that 
		\begin{align} 
		\frac{1}{2} \frac{|\vm^{\underline{\lambda}}|^2}{\vr^{\underline{\lambda}}}
		(t, \cdot) &\leq \chi(t) \ \mbox{for a.a.}\ t > \tau,\br
		\frac{1}{2} \frac{|\vm^{\Ov{\lambda}}|^2}{\vr^{\Ov{\lambda}}}
		(t, \cdot) &\geq \chi(t) \ \mbox{for a.a.}\ t > \tau,\ \frac{1}{2} \frac{|\vm^{\Ov{\lambda}}|^2}{\vr^{\Ov{\lambda}}}
		(t, \cdot) \to 0 \ \mbox{as}\ t \to \infty.		
		\label{C2}
		\end{align}	
		
	\end{itemize}
\end{Theorem}	

Property \eqref{C2} means that the solution approaches asymptotically the equilibrium state 
\[
(\vr_s, \vm_s , s_s),\ 
\vr_s = \vr_0,\ \vm_s = 0, \ p(\vr_0, s_s) = \frac{2}{3} \Lambda.
\]
The rate of convergence can be arbitrarily fast/slow depending on the choice of 
the profile $\chi$.

\begin{Remark} \label{final}
	
As matter of fact, there are solutions stabilizing to equilibrium in an arbitrarily short positive time $\tau$. Indeed it is enough to extend 
the solutions claimed in Proposition \ref{rP2} as 
\[
\vr (t, \cdot) = \vr_0,\ \vm(t, \cdot) = 0, \ s(t, \cdot) = \widehat{s}(t, \cdot)\ \mbox{for}\ t \geq \tau, 
\]
where $\widehat{s}(t, \cdot) \geq s(\tau- , \cdot)$ is chosen so that the total energy is conserved, meaning 
\[
\vr_{0,n} e(\vr_{0,n}, \widehat{s}) = 
\vr_{0,n} e(\vr_{0,n}, s_{0,n}) + \Lambda - \frac{d}{2} p(\vr_{0,n}, s_{0,n}) 
= \Lambda > \vr_{0,n} e(\vr_{0,n}, s_{0,n})\ \mbox{in}\ \Omega_n, 
\]
cf. \eqref{r21}.

\end{Remark}

\section{Acknowledgments}

The work of E.F.~was partially supported by the Czech Sciences Foundation (GA\v CR), Grant Agreement 24--11034S. The Institute of Mathematics of the Czech Academy of Sciences is supported by RVO:67985840. E.F. is a member of the Ne\v cas Center for Mathematical Modelling. A.A. is member of GNFM - INdAM.

The paper was elaborated during the visit of E.F. to Dipartimento di Matematica e Fisica of the University of Campania ``Vanvitelli''. 

\def\cprime{$'$} \def\ocirc#1{\ifmmode\setbox0=\hbox{$#1$}\dimen0=\ht0
	\advance\dimen0 by1pt\rlap{\hbox to\wd0{\hss\raise\dimen0
			\hbox{\hskip.2em$\scriptscriptstyle\circ$}\hss}}#1\else {\accent"17 #1}\fi}


\end{document}